\documentclass[12pt,oneside,english]{amsart}
\usepackage[T1]{fontenc}
\usepackage[utf8]{inputenc}
\usepackage{geometry}
\usepackage{fancyhdr}
\usepackage{color}
\usepackage[foot]{amsaddr}
\usepackage{amstext}
\usepackage{amsthm}
\usepackage{amssymb}
\usepackage{mathdots}
\usepackage{esint}
\usepackage{listings}
\usepackage{setspace}
\usepackage[noadjust]{cite}
\usepackage{subcaption}
\usepackage[ruled]{algorithm2e}

\makeatletter
\numberwithin{equation}{section}
\numberwithin{figure}{section}
\theoremstyle{plain}
\newtheorem{thm}{\protect\theoremname}
  \theoremstyle{plain}
  \newtheorem{lem}[thm]{\protect\lemmaname}
  \theoremstyle{remark}
  \newtheorem*{rem*}{\protect\remarkname}
  \theoremstyle{remark}
  \newtheorem{rem}[thm]{\protect\remarkname}
  \theoremstyle{plain}

\usepackage[english]{babel}
\usepackage{fancyhdr}% http://ctan.org/pkg/fancyhdr
\definecolor{Code}{rgb}{0,0,0}
\definecolor{Decorators}{rgb}{0.5,0.5,0.5}
\definecolor{Numbers}{rgb}{0.5,0,0}
\definecolor{MatchingBrackets}{rgb}{0.25,0.5,0.5}
\definecolor{Keywords}{rgb}{0,0,1}
\definecolor{self}{rgb}{0,0,0}
\definecolor{Strings}{rgb}{0,0.63,0}
\definecolor{Comments}{rgb}{0,0.63,1}
\definecolor{Backquotes}{rgb}{0,0,0}
\definecolor{Classname}{rgb}{0,0,0}
\definecolor{FunctionName}{rgb}{0,0,0}
\definecolor{Operators}{rgb}{0,0,0}
\definecolor{Background}{rgb}{0.98,0.98,0.98}

\providecommand{\lemmaname}{Lemma}

\providecommand{\theoremname}{Theorem}
\providecommand{\corollaryname}{Corollary}
\newcommand{\abs}[1]{\ensuremath{|#1|}}

\newcommand{\norm}[2]{\ensuremath{|\!|#1|\!|_{#2}}}

\newcommand{\Norm}[2]{\ensuremath{\left|\!\left|#1\right|\!\right|_{#2}}}

\newcommand{\braket}[2]{\langle #1 | #2 \rangle}

\renewcommand{\d}[1]{\ensuremath{\textnormal{d}#1}}

\newcommand{\cF}{\mathcal{F}}
\newcommand{\cG}{\mathcal{G}}
\newcommand{\cH}{\mathcal{H}}

\newcommand{\cK}{\mathcal{K}}
\newcommand{\cL}{\mathcal{L}}

\newcommand{\cN}{\mathcal{N}}

\newcommand{\cP}{\mathcal{P}}

\newcommand{\cR}{\mathcal{R}}
\newcommand{\cS}{\mathcal{S}}

\usepackage{babel}
\providecommand{\lemmaname}{Lemma}
  
  \providecommand{\remarkname}{Remark}
\providecommand{\theoremname}{Theorem}

\usepackage{babel}
\providecommand{\lemmaname}{Lemma}
  
\providecommand{\theoremname}{Theorem}

\usepackage{babel}
\providecommand{\theoremname}{Theorem}

\usepackage{babel}
\providecommand{\theoremname}{Theorem}

\usepackage{mathtools}

\DeclarePairedDelimiterX{\inp}[2]{\langle}{\rangle}{#1, #2}

\makeatother

\theoremstyle{plain}

  \theoremstyle{plain}
    
  \theoremstyle{plain}
      \newtheorem{example}[thm]{\protect\examplename}
  \theoremstyle{plain}

\providecommand{\lemmaname}{Lemma}

\providecommand{\theoremname}{Theorem}
\providecommand{\definitionname}{Definition}
\providecommand{\remarkname}{Remark}
\providecommand{\examplename}{Example}
\DeclareMathOperator*{\argmin}{arg\,min}

\begin{document}
\sloppy

\title{An exact kernel framework for spatio-temporal dynamics}

\author{Oleg Szehr}
\thanks{The author expresses his thanks to an anonymous referee who pointed out the IDW approach to spatio-temporal regression.}
\author{Dario Azzimonti}
\author{Laura Azzimonti}
\address{Dalle Molle Institute for Artificial Intelligence (IDSIA) - SUPSI/USI, Manno, Switzerland}

\email{oleg.szehr@idsia.ch}
\email{dario.azzimonti@idsia.ch}
\email{laura.azzimonti@idsia.ch}


\begin{abstract}
A kernel-based framework for spatio-temporal data analysis is introduced that applies in situations when the underlying system dynamics are governed by a dynamic equation. The key ingredient is a representer theorem that involves time-dependent kernels. Such kernels occur commonly in the expansion of solutions of partial differential equations. The representer theorem is applied to find among all solutions of a dynamic equation the one that minimizes the error with given spatio-temporal samples. This is motivated by the fact that very often a differential equation is given a priori (e.g.~by the laws of physics) and a practitioner seeks the best solution that is compatible with her noisy measurements. Our guiding example is the Fokker-Planck equation, which describes the evolution of density in stochastic diffusion processes. A regression and density estimation framework is introduced for spatio-temporal modeling under Fokker-Planck dynamics with initial and boundary conditions.
\end{abstract}

\maketitle

\section{Introduction}
Spatio-temporal processes occur throughout the scientific disciplines with wide-spread applications in areas such as physics, economy, climate, ecology,...~\cite{Wikle2011,Wikle2019}. The respective data analysis requires the modeling of spatial and temporal dynamics. For processes that occur in nature, the spatio-temporal interaction is typically governed by a dynamic equation. In many cases the underlying dynamics are known a priori (e.g.~from the laws of physics or economy) and a practitioner/ an algorithm searches for a solution that optimally matches spatio-temporal data. Learning exact solutions of dynamic equations is thus a key goal of spatio-temporal modeling. This article introduces a new kernel-based learning theory that applies to spatio-temporal data analysis under dynamic equation constraints. The constraints are  implemented in terms of positive definite kernels that solve the dynamic equation~\emph{exactly}. This confers the main advantages of our method~\emph{1)} coverage of domain constraints in terms of initial and boundary conditions and~\emph{2)} exact modeling of the spatial and temporal dynamics of a process. Positive definite kernels play a fundamental role as building blocks of kernel-regression and density estimation in non-parametric statistics~\cite{Silverman1986} and of support vector machines~\cite{SS2002} in machine learning. The representer theorem can be seen as the cardinal result in the application of kernel-based methods. It states that the minimizer of an empirical risk functional can be expanded in terms of kernels evaluated at data samples alone. On the technical side the main innovation behind our framework lies in the introduction of a dynamic risk-minimization framework and a respective representer theorem that applies in the spatio-temporal setting.

%This article introduces a new kernel-based learning theory that applies to spatio-temporal data analysis under dynamic equation constraints with a risk-minimization framework and a generalized, time-dependent representer theorem. As compared to the classical representer theorem, this setup involves a $t$-family of reproducing kernel Hilbert spaces (RHKS), leaving open the question how to combine information from different spaces.

Formally, the present article investigates the learning of a time-dependent process $f=f(x,t)$ from samples $\cS^{(k)}$ observed at different times $t_k$, $k\in\{1,...,T\}$.
The main underlying assumption is that the function $f$ can be expanded in terms of a convergent series of time-\textbf{de}pendent positive kernels
\begin{equation}
\label{main}f(x,t)=\sum_{\mu=1}^\infty a_\mu\cK_t(x,x_\mu),
\end{equation}
with time-\textbf{in}dependent coefficient $\{a_\mu\}_{\mu=1}^\infty$. Functions of the form~\eqref{main} occur commonly as solutions of time-dependent partial differential equations (PDEs). Two flavors of the learning problem are studied here: 
\begin{enumerate}
\item \emph{Spatio-temporal Kernel Regression:} Suppose at times $t_1,...,t_T$ noisy data samples $\cS^{(k)}=\{(x_i^{(k)},y_i^{(k)})\}_{i=1,...,N}$ are taken from a function that evolves according to~\eqref{main}. What is the optimal fit  to the total sample $\cS=\cup_{k=1}^T\cS^{(k)}$ among the space of solutions?
\item \emph{Spatio-temporal Kernel Density estimation:} Samples $\cS^{(k)}=\{x_i^{(k)}\}_{i=1,...,N}$ are taken from an evolving density at times $t_1,...,t_T$. In this case~\eqref{main} reflects a distributional embedding~\cite{Smola2007}, see below. We construct a kernel estimator for the evolving density that takes account of all samples simultaneously.
\end{enumerate}
Our guiding application is the learning of solutions of time-dependent diffusion equations of the form
\begin{equation}
L u(x,t) = \partial_t u(x,t),\label{FokkerPlanck}
\end{equation}
where $L$ is a (self-adjoint) operator of Fokker-Planck type, see below for details. $u$ represents either a regression function or density. This is motivated by the ubiquitous occurrence of It\^{o} diffusion processes throughout scientific disciplines, see e.g.~\cite{Allen2010,Risken1989,Brigo2001}. The simplest instance of~\eqref{FokkerPlanck} is the heat equation
\begin{equation}
\partial_{xx} u(x,t) = \partial_t u(x,t)\label{heatEquation},
\end{equation}
which describes the density of the ordinary Brownian motion.
Its solutions can be written in the form~\eqref{main} with Gaussian kernels 
\begin{equation*}
\cK^{Gauss}_t(x,x_i)=\frac{1}{\sqrt{2\pi t}}e^{-\frac{(x-x_i)^2}{2t}}.
\end{equation*}
Other common examples include the geometric Brownian motion as a model for stock prices and the Hull-White model (Ornstein-Uhlenbeck process) for interest rates in economy. 

\begin{figure}
%\begin{figure}[H]
\centering
\begin{subfigure}[b]{0.25\textwidth}
\includegraphics[width=\linewidth]{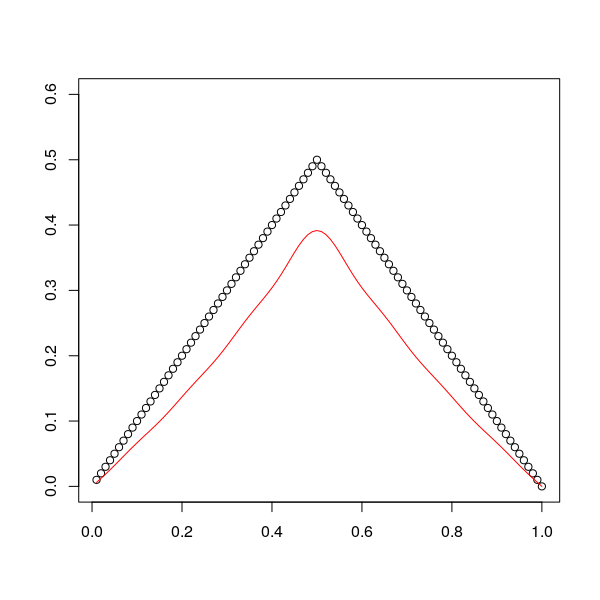}
\caption{$t_1=0.01$}
\end{subfigure}
\begin{subfigure}[b]{0.25\textwidth}
\includegraphics[width=\linewidth]{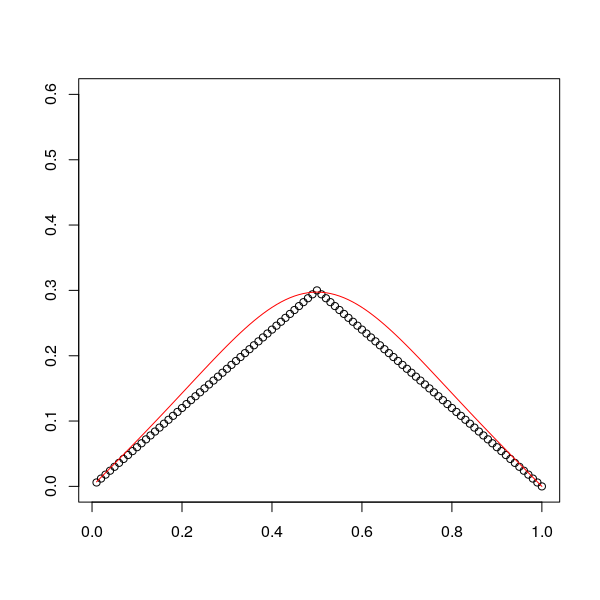}
\caption{$t_2=0.02$}
\end{subfigure}
\begin{subfigure}[b]{0.25\textwidth}
\includegraphics[width=\linewidth]{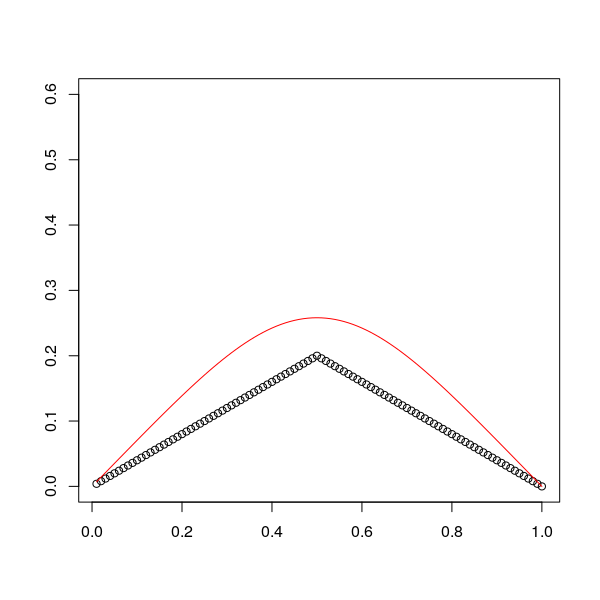}
\caption{$t_2=0.03$}
\end{subfigure}
\caption{Optimal match to temperature measurements by $100$ sensors distributed at equal distances over a metal rod. Three series of measurements are taken at times $t_1, t_2, t_3$. It is assumed that temperature, $u$, evolves according to~\eqref{heatEquation} and $u=0$ at the boundaries of the metal rod (Dirichlet boundary condition). Black circles represent exemplary measurement data. Red lines depict a regression function that minimizes $l^2$-risk.}\label{introExample}
\end{figure}

In the context of PDEs, initial and boundary conditions play a crucial role. For regression they can lead to {over-determined} learning problems, i.e.~no solution of the PDE complies simultaneously with boundary conditions and measured data. In practice this can occur as a result of noise. In such cases our method can be applied to balance between mismatches in boundary conditions and measurement errors. For illustration Figure~\ref{introExample} shows an example of an over-determined spatio-temporal regression with solutions of the heat equation.

\begin{figure}
\centering
\begin{subfigure}[b]{0.24\textwidth}
\includegraphics[width=\linewidth]{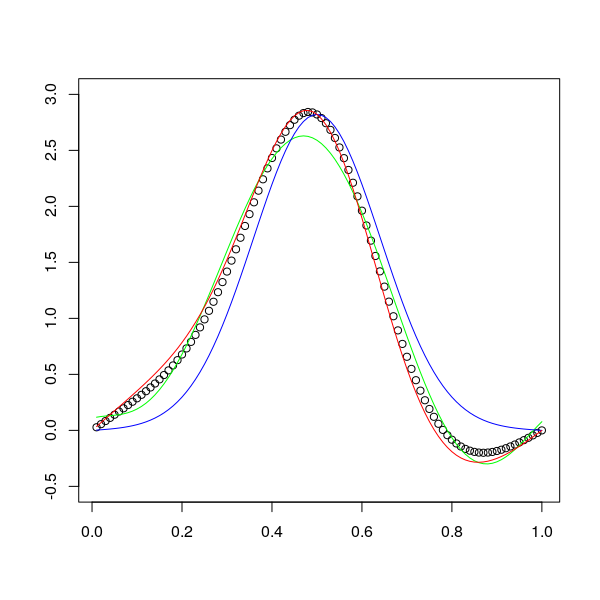}
\caption{$t_1=0.02$}
\end{subfigure}
\begin{subfigure}[b]{0.24\textwidth}
\includegraphics[width=\linewidth]{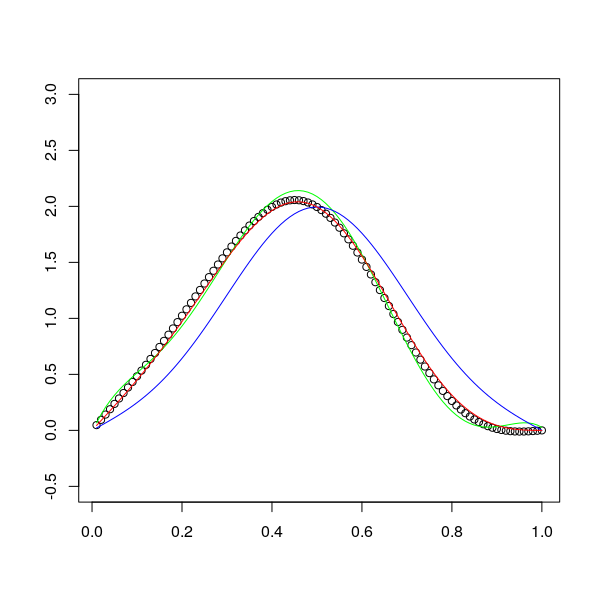}
\caption{$t_2=0.03$}
\end{subfigure}
\begin{subfigure}[b]{0.24\textwidth}
\includegraphics[width=\linewidth]{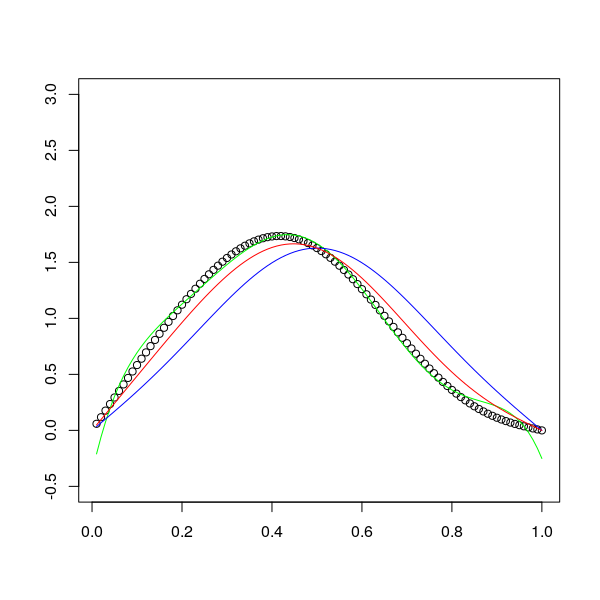}
\caption{$t_3=0.03$}
\end{subfigure}
\begin{subfigure}[b]{0.24\textwidth}
\includegraphics[width=\linewidth]{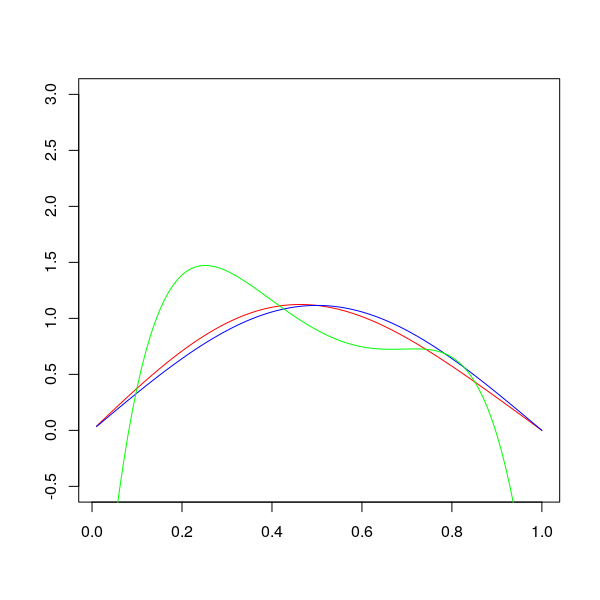}
\caption{$t_4=0.06$}
\end{subfigure}
\caption{Comparison of PDE and IDW approaches under Dirichlet boundary conditions. Circles indicate exemplary noisy data. A regression model is calibrated to data at times $t_1,t_2,t_3$, see (A)--(C). (D) shows respective predictions for time $t_4$. Blue lines indicate the exact evolution of temperature. Green lines are obtained from a Gaussian kernel predictor with bandwidth $0.45$. Red lines are obtained from our PDE-based prediction.}\label{figureHeatIDW}
\end{figure}

\textbf{Background on spatio-temporal models:} A common method to account for the temporal discrepancy of data is~{inverse distance weighting} (IDW)~\cite{Wikle2011,Wikle2019}. The main idea is that values at unknown points in time are computed as weighted averages of the available data snap-shots, where higher weight is given to temporally closer data\footnote{This is known as 'Shepard's method' in spatial data analysis. Another name, 'Tobler's law', is common also for spatio-temporal problems~\cite{Miller2004} and stems from geography.}. This form of IDW can be viewed as a spatio-temporal kernel method, where the weights are given by
$$w_{(ik)(jl)} = \cK((x_i^{(k)},t_k),(x_j^{(l)},t_l)).$$
Thus a 'straight-forward' approach to spatio-temporal modeling is to follow the IDW paradigm by incorporating time-coordinates as part of the observations into a given kernel. In reality, however, the evolution of a system is much less 'symmetric' with respect to space and time and usually governed by a dynamic equation. This article introduces spatio-temporal modeling, where a time-dependent kernel $\cK_t$ satisfies a PDE exactly, as an alternative to the IDW approach. Finally, hierarchical dynamic spatio-temporal models~\cite{Wikle2010} provide a framework that can capture a much wider variety of non-separable models than kernel IDW but our kernels exactly follow the systems dynamic equation. The key advantages of our method are~\emph{1)} that domain constraints are covered naturally in terms of initial and boundary conditions and~\emph{2)} that it allows for an exact modeling of the spatial and temporal dynamics of a process. This is of particular interest for prediction, e.g.~in the presence of data drifts that are captured within Fokker-Planck dynamics. Diffusion-based kernel density estimation~\cite{Botev2010} employs the Fokker-Planck equation as a resource for the construction of kernels under domain constrains and pilot density estimates. It appears natural to study such kernels also in the context of spatio-temporal modeling. Hence, the first main contribution of this article lies in the observation that spatio-temporal inter-dependence in many cases is much better captured via Fokker-Planck dynamics with time-dependent kernel than a fixed multivariate kernel with optimized bandwidth. Figure~\ref{figureHeatIDW} shows an example of spatio-temporal prediction comparing IDW and PDE approaches, see Example~\ref{exampleIDW} for details. Remark the key feature of the new PDE predictor, that high error levels in training data, (A)--(C), are 'smoothed' out to provide an accurate prediction (D) that is compliant with the underlying heat equation.

\textbf{Background on PDEs:} Kernel expansion techniques are commonly employed as a \emph{meshless} (as opposed to \emph{finite difference}) method for the numerical approximation of PDE solutions~\cite{Belytschko1996}. For time-{in}dependent PDEs, kernel discretization represents the trial function entirely as a sum of 'data-nodes'. This approach follows a classical line of~\cite{Trefftz1926} to use trial functions that satisfy a PDE exactly. For time-dependent PDEs, meshless kernel-based methods usually rely on a time-independent spatial sum representation, but the coefficients explicitly depend on time. In contrast the article at hand expands solutions as a series of time-dependent kernels $\cK_t$ that solve~\eqref{main} exactly (following more rigidly the lines of~\cite{Trefftz1926}). While this adds complexity to the structure of kernels it conveniently eliminates the time integration from solving the PDE and allows for a natural representer theorem. For an introduction to kernel methods with focus on meshless PDE solutions and machine learning see~\cite{Schaback2006}. Our work is inspired by the recent article~\cite{Hon2015}, which studies the numerical approximation of solutions of~\eqref{FokkerPlanck} when $L$ is a spatial, elliptic operator. The important representation of solutions in the form~\eqref{main} is introduced there as a method for meshless approximation. While~\cite{Hon2015} focuses on kernel-based interpolation of initial conditions, here, the technique is extended to learn from spatio-temporal data. 
%Our focus the Fokker-Planck equation~\eqref{FokkerPlanck} following well-established spectral methods as described in~\cite{Risken1989}. 
Thus the second main contribution of the present article lies in the recognition that the representation~\eqref{main} is particularly suited for statistics and machine-learning applications with spatio-temporal samples and in the derivation of a dynamic representer theorem for spatio-temporal models.

%We illustrate the practical relevance of the introduced framework by the following real-world examples:
%
%\begin{enumerate}
%\item \emph{Regression:} Temperature is measured periodically by sensors at different spatial locations. The objective is to predict temperature at new spaces and times. Boundary conditions are given e.g.~by the shape of the domain of the temperature distribution.
%\item \emph{Density estimation:} A diffusion process of It\^{o} type describes the price of a financial asset. Samples are taken on daily basis, given available price information. The structure of the It\^{o} process mirrors some micro-economic laws, e.g.~no-arbitrage requirements, etc...
%A common boundary condition is that prices shall be non-negative.
%\end{enumerate}
%
%
%
%
%
\subsection{Recap on statistical learning}\label{recapStatisticalLearning}
We recapitulate some concepts of supervised learning. For details we refer the reader to the survey article~\cite{Luxburg2008} and references therein.
Let $\cF\subset\mathbb{Y}^\mathbb{X}$ be a set of functions $f:\mathbb{X}\rightarrow\mathbb{Y}$.
Given a fixed but unknown probability distribution $\mathbb{P}$ on $\mathbb{X}\times\mathbb{Y}$, the risk of $f\in\cF$ is
\begin{equation*}
R(f):=\int_{\mathbb{X}\times\mathbb{Y}}loss(f(x),y)\d\mathbb{P}(x,y).
\end{equation*}
$loss(y_1,y_2)$ measures the discrepancy between $y_1,y_2\in\mathbb{Y}$ and for simplicity is assumed to be convex in what follows. Risk measures by how much (accumulated over $\mathbb{X}$) the learned responses $f(x)$ deviate from the ground truth $y$. The goal of learning is to identify $f^*$ that has as small as possible risk within the constraints imposed by $\cF$,
\begin{equation*}
f^*=\argmin_{f\in\cF} R(f).
\end{equation*}
In practice, it is impossible to determine $f^*$ at least for the following reasons:

\emph{(1)} Without structural constraints on $\cF$, the minimization is intractable.

\emph{(2)} $\mathbb{P}$ is not known.

As a consequence of ${(1)}$ an a priori structure is commonly imposed on $\cF$. Regarding point ${(2)}$, a data sample $\cS=\{(x_i,y_i)\}_{i=1,...,N}\in(\mathbb{X}\times\mathbb{Y})^N$ is taken i.i.d.~from $\mathbb{P}$. Given the sample, a natural approach is to minimizes the empirical risk
\begin{align}
R_{emp}^{\cS}(f):=\frac{1}{N}\sum_{i=1}^N loss(f(x_i),y_i).\label{empiricalRisk}
\end{align}

The classical representer theorem asserts that if this minimization is performed over a reproducing kernel Hilbert space (RKHS), then any optimal function is a \emph{finite} sum of kernels evaluated at sampled data points~\cite{Wahba1999,Christmann2008}. This reduces the infinite-dimensional optimization problem to an $N$-dimensional sub-space making it accessible to numerical optimization. For completeness, the classical representer theorem asserts the following.

\begin{thm}[Classical representer theorem]\label{repThmClassical}
Suppose that $\cK$ is a symmetric, positive definite kernel on a set $\mathbb{X}\neq\emptyset$. Let
$\cH=\left\{f=f(x)\ |\ f(x)=\sum_{\mu=1}^\infty a_\mu \cK(x,x_\mu)\right\}$ denote the RKHS spanned by convergent sums of kernels.
Let $\cS=\{(x_i,y_i)\}_{i=1,...,N}$ denote a data sample taken from $(\mathbb{X}\times\mathbb{Y})^N$. Then
%\cup_{\mu=1}^\infty\{x_\mu\}$. Then
%
%
\begin{align*}
\inf_{f\in\cH}R_{emp}^\mathcal{\cS}(f)=\min_{f\in\cH_{rep}}R_{emp}^\mathcal{\cS}(f),
\end{align*}
where
$$\cH_{rep}=\left\{f=f(x)\ |\ f(x)=\sum_{\mu=1}^N a_{\mu} \cK(x,x_\mu)\right\}.$$
\end{thm}
\begin{rem}
In what follows we assume, without saying explicitly, that series of the form $f(x)=\sum_{\mu=1}^\infty a_\mu \cK(x,x_\mu)$ are convergent in the norm of the Hilbert space. This is justified by the fact that we only consider $\cK$ within its \emph{native} (i.e.~Moore-Aronszajn) RKHS~\cite{Schaback2006}. Similarly series involving time-dependent kernels $\cK_t$ converge for fixed $t$ within the respective RKHS.  
\end{rem}
The article at hand introduces a time parameter to risk minimization. This involves a generalization of the empirical risk functional~\eqref{empiricalRisk} and an extension of Thm.~\ref{repThmClassical} beyond the setting of a single RKHS. 

For completeness we mention the important topic of regularization. To prevent overfitting, a 'penalty term' (regularizer) is added to the empirical risk functional. This has been introduced in~\cite{Kimeldorf1971} and extended to regularizers of the form $g(\norm{f}{})$ with strictly increasing $g$ in~\cite{Schoelkopf2001}. The article \cite{Argyriou2009} investigates a matrix-valued setup and derives conditions for a respective representer theorem. In the interest of space we do not cover regularization here, although regularized versions of our main theorems~\ref{thmRep},~\ref{repThmDEq} and~\ref{repFullDensityLearning} clearly hold.

\subsection{Recap on the Fokker-Planck equation}\label{sec:recapFP}
We recapitulate some standard results from the theory of the Fokker-Planck equation. Details and derivations can be found in the monograph~\cite{Risken1989}.
The Fokker-Planck equation plays a central role in the theory of stochastic processes: a normalized solution of the equation represents the probability density of a random motion that follows a stochastic process of It\^{o} type. The simplest physical setting is that of a particle undergoing a diffusion process. The simplest economic setting is that of a random market price of a financial asset. 
To set the stage, suppose a one-dimensional It\^{o} process
\begin{equation*}
\d X_{t}=\mu (X_{t})\,dt+\sigma (X_{t})\,dW_{t}
\end{equation*}
with time-independent\footnote{It\^{o} processes with time-dependent drift or diffusion require a more general discussion although some of our findings still apply.} drift $\mu=\mu(x)$ and diffusion $D=D(x)=\sigma^{2}(x)/2$. The Fokker-Planck equation for the probability density ${\rho(x,t)}$ of the random variable $X_{t}$ is
\begin{equation*}{\partial_t\rho(x,t)=-\partial_x\left[\mu\rho(x,t)\right]+\partial_{xx}\left[D\rho(x,t)\right]}.
\end{equation*}
For brevity it is common to introduce the Fokker-Planck operator
$$L:=-\partial_x\mu+\partial_{xx}D$$
writing the differential equation in the form~\eqref{FokkerPlanck}. An important role in the study of boundary conditions of this PDE is played by the occurrence of a \emph{stationary density} $\rho_{fix}(x)$, i.e.
$$\lim_{t\to\infty}\rho(x,t)=\rho_{fix}(x).$$
If a stationary solution exists then it can be written in the form
$$\rho_{fix}(x)=Ce^{-\Phi(x)},$$
with a potential of the form $\Phi(x)=\ln(D(x))-\int^x\frac{\mu(x')}{D(x')}\d x'$ and a normalization constant $C$. The potential $\Phi$ is a natural point to incorporate boundary conditions of Dirichlet type: If the potential jumps to an infinite value at $x\geq x_0$ then the density is confined to an interval $x<x_0$. More generally, if the potential satisfies appropriate growth conditions then there is a unique stationary density of the above form. In this situation it is common to introduce a scalar product
$$\braket{f}{g}_{e^\Phi}=\int f(x){g}(x)e^{\Phi(x)}\d x.$$ The respective Hilbert function space is $L^2(e^{\Phi(x)}\d x)$. A quick computation (that is detailed out in~\cite[Section~5.4]{Risken1989}) shows that $L$
is self-adjoint with respect to $\braket{\cdot}{\cdot}_{e^{\Phi(x)}}$ under various types of boundary conditions. Hence, our working assumption is:

\emph{(A) The operator $L$ is a negative\footnote{The assumption of negativity can be relaxed. If positive eigenvalues are present, a time horizon $T$ is fixed for solutions of the PDE. The requirement is that the eigenfunction expansion of the kernel~\eqref{eigenExpansion} is convergent at $T$.} self-adjoint operator on the Hilbert space $\cH=L^2(e^{\Phi(x)}\d x)$ for a certain potential $\Phi$.}

In what follows only \emph{(A)} itself is relevant rather than its exact origin. Conditions that are sufficient for \emph{(A)} are described in the literature, see~\cite[Section~5.4]{Risken1989} for an introductory treatment. The spectrum of $L$ can be discrete or continuous or both, a spectral decomposition follows from the spectral theorem. For convenience we focus on the discrete spectrum, but the discussion can be generalized using standard tools of functional analysis. In this situation, a complete orthonormal basis $\{\varphi_n\}_n$ (wrt.~$\braket{\cdot}{\cdot}_{e^{\Phi(x)}}$) diagonalizes $L$, 
$$L\varphi_n=-\lambda_n\varphi_n,\quad \lambda_n\geq0, \quad n=1,2,...$$
We introduce the kernel $\cK^{FP}_t(x,x')$ by the formula 
\begin{align}
\cK^{FP}_t(x,x')= \sum_{n=1}^\infty e^{-\lambda_nt}\varphi_n(x)\varphi^{T}_n(x')\label{eigenExpansion}
\end{align}
and notice that it solves~\eqref{FokkerPlanck} exactly. Notice also that the completeness relation for the eigenfunctions reads
$\delta(x-x')=\sum_n\varphi_n(x)\varphi^T_n(x')$, i.e.~the kernel $\cK^{FP}_t(x,x')$ is nothing but the Green's function of the Fokker-Planck equation. One can obtain the general solution of the PDE with initial conditions $f(x,0) = g(x)$ as a convolution integral 
\begin{align} f(x,t)=\int_{\mathbb{X}}g(x') \cK^{FP}_t (x,x')\d x'.\label{convIntegral}
\end{align}
The representation for $\cK^{FP}_t(x,x')$ is of the same form as discussed in~\cite{Hon2015} for the example of the heat equation.

\begin{lem}For any $t>0$ $\cK^{FP}_t(x,y)$ is a symmetric, positive-definite kernel on $\mathbb{X}$ and for any fixed $t$ the kernel $\cK^{FP}_t$ generates a unique RKHS $\cH_{\cK_t^{FP}}$.\label{positiveKernel}
\end{lem}
The simple proof is shown in the appendix.
\begin{example}(Ornstein-Uhlenbeck process) To illustrate the framework, consider the Ornstein-Ulenbeck (OU) process, which plays an important role in interest rate models in financial Mathematics~\cite{Brigo2001}. For fixed $\theta,\sigma>0$, the zero-mean OU process is characterized by the stochastic differential equation
$$\d X_t = -\theta X_t\d t+\sigma\d W_t.$$
The process is mean-reverting in the sense that if $X_t<0$ then the drift term $-\theta X_t$ is positive and if $X_t>0$ then the drift term is negative. $\theta$ determines the speed of mean reversion. The associated Fokker-Planck equation is
$${\partial_t\rho=\theta\partial_x\left(x\rho\right)+\frac{\sigma^2}{2}\partial_{xx}\rho},$$
i.e.~$L=\theta\partial_x(x\cdot)+\frac{\sigma^2}{2}\partial_{xx}(\cdot)$.
The eigenvalues are $\lambda_n=n$, $n\geq0$ with eigenfunctions
$$\varphi_n(x)=\frac{1}{\sqrt{\sqrt{2\pi}n!}}e^{-\frac{1}{2}\left(x/\sqrt{\frac{\sigma^2}{2\theta}}\right)^2}H_n\left(x/\sqrt{\frac{\sigma^2}{2\theta}}\right).$$ Here $H_n(x)=(-1)^ne^{x^2/2}\frac{d^n}{dx^n}e^{-x^2/2}$ denote the Hermite polynomials. The stationary density corresponds to the eigenvalue $n=0$, that is~$\rho_{fix}=\cN\left(0,\frac{\sigma^2}{2\theta}\right).$
\end{example}

\subsection{Recap on kernel density estimation}\label{recapDensityEstimation}
Given $N$ independent realizations $X_1,...,X_N$ from an unknown continuous probability density function $\rho$ on $\mathbb{X}$ and a positive kernel $\cK$, the kernel density estimator is
\begin{align}
\hat{\rho}(x)=\frac{1}{N}\sum_{i=1}^N\cK(x,X_i).\label{kernelDensityEstimator}
\end{align}
The methodology developed here relies on the \emph{kernel embedding of distributions} of~\cite{Smola2007}. Given a kernel $\cK$ on $\mathbb{X}\times\mathbb{X}$ the distributional embedding $\varepsilon$ of a density $\rho$ into the RKHS $\cH_\cK$ of $\cK$ is
\begin{align*}
[\varepsilon(\rho)](x)=\int_{\mathbb{X}} \cK(x,x')\rho(x')\d x'\in\cH_\cK.
\end{align*}

On the technical side it should be mentioned that the integral is assumed to exist and that the embedding is assumed to be injective. A detailed discussion of the technical backbone can be found in~\cite{Smola2007, Steinwart2002}. The embedding allows one to operate on distributions using Hilbert space concepts such as inner products and linear projections. For us it serves two purposes.

\emph{1)} It is the entrance point to interpret density estimation within the risk minimization framework.

\emph{2)} In the context of evolving densities and time-dependent kernels, the embedding naturally mirrors the structural properties of the PDE that governs the evolution. 

The kernel mean estimator (KME) is the minimizer of the empirical risk functional (compare~\eqref{empiricalRisk}) of the embedding~\cite{Muandet2014, STC2004},
\begin{align}\label{KME}
\cR_{KME}(\varepsilon)=\frac{1}{N}\sum_{i=1}^N \Norm{\varepsilon-\cK(\cdot,x_i)}{RKHS}^2,
\end{align}
where $\Norm{\cdot}{RKHS}$ denotes the norm of the RKHS. The representer theorem, Thm.~\ref{repThmClassical}, applies and yields a representation of the KME in the form of an ordinary kernel density estimator
$\hat{\rho}(x)=\sum_{i=1}^Na_i \cK(x,X_i)$. In the absence of a regularizer, equal weights are optimal~\cite[Prop.~5.2]{STC2004}, which yields~\eqref{kernelDensityEstimator}.
%
%
%
%

%In this article we study density estimators that satisfy the Fokker-Planck equation. We develop the methodology to incorporate samples taken at times $\{t_k\}_{k=1,...,T}$ from an evolving density simultaneously. For the Fokker-Planck kernels $\cK^{FP}$ our method fits naturally with the representation of the solutions of the Fokker-Planck equation in terms of a convolution integral. %We interpret $\rho(x)$ as an initial density for a stochastic process and we are looking for an estimator $\hat\rho$ for $\rho$.
%
%\begin{example}
%
%
%
%
%
%
%Our guiding example is the Gaussian kernel density estimator
%
%
%\begin{align*}
%\hat{\rho}(x,t)=\frac{1}{N}\sum_{i=1}^N\cK_t^{Gauss}(x,X_i),
%\end{align*}
%
%
%
%which is a sum of Green's functions of the heat equation. Notice that this estimator is intrinsically time-dependent and solves the Fokker-Planck equation exactly. For each fixed $t_k$ this is nothing but the ordinary kernel density estimator~\eqref{kernelDensityEstimator}.
%\end{example}

%
%\subsection{Outline}
%In Section~\ref{supervisedLearning} we discuss kernel-regression. The space of admissible functions $\cF$ is the space of solution to the Fokker-Planck equation with boundary conditions.
%
%Section~\ref{densityLearning} extends the developed framework to density learning. In this situation a density that evolves according to the Fokker-Planck equation is estimated from samples that are taken at different times.

\section{Results}

\subsection{Supervised Learning of time-dependent functions}\label{supervisedLearning}

Let noisy data samples $\cS^{(k)}=\{(x_i^{(k)},y_i^{(k)})\}_{i=1,...,N}$ be taken at different times $t_k$, $k\in\{1,...,T\}$ and let $\cS=\cup_{k=1}^T\cS^{(k)}$ denote the full data sample over all times. Let $f=f(x,t)$ be a time-dependent function to be estimated via kernel regression from a set of admissible functions $\cG$, using information contained in $\cS$. We define the empirical risk of $f$ over $\cS$ as an equally weighted sum of the empirical risks of $f(\cdot,t_k)$ over $\cS^{(k)}$,
\begin{align}
R_{emp}^\mathcal{\cS}(f)&=\frac{1}{T}\sum_{k=1}^TR_{emp}^{\mathcal{\cS}_k}(f(\cdot,t_k)).\label{regressionLearning}
\end{align}
This generalizes~\eqref{empiricalRisk} to the case that $T>1$ but can also be viewed as a special case in the sense that $\cup_{i,k=1}\{((x_i^{(k)},t_k),y_i^{(k)})\}$ is interpreted as the underlying data sample. We begin by studying the situation when samples are taken from a countable set $\mathbb{X}$ assuming that $\cG$ is composed of functions of the form $f(x,t)=\sum_{\mu=1}^\infty a_\mu \cK_t(x,x_\mu)$,
with symmetric and positive kernels $\cK_t$. 
%The underlying assumption is that equal weights are given to samples obtained at different times $t_k$.

%
%
\begin{thm}[Time-dependent representer theorem]\label{thmRep}
Suppose that $\cK_t$ is a symmetric, positive definite kernel on a set $\mathbb{X}$ for every $t>0$. Let
$$\cG=\left\{f=f(x,t)\ |\ f(x,t)=\sum_{\mu=1}^\infty a_\mu \cK_t(x,x_\mu)\right\}$$ and suppose that at times $t_k$, $k\in\{1,...,T\}$, data samples $\cS^{(k)}=\{(x_i^{(k)},y_i^{(k)})\}_{i=1,...,N}$ are measured and that $\cup_{t,k}\{x_i^{(k)}\}\subset\mathbb{X}$. Then
%\cup_{\mu=1}^\infty\{x_\mu\}$. Then
%
%
\begin{align*}
\inf_{f\in\cG}R_{emp}^\mathcal{\cS}(f)=\min_{f\in\cG_{rep}}R_{emp}^\mathcal{\cS}(f),
\end{align*}
where
\begin{small}
$$\cG_{rep}=\left\{f=f(x,t)\ |\ f(x,t)=\sum_{k=1}^T\sum_{i=1}^N a_{i,k} \cK_t(x,x_i^{(k)})\right\}.$$
\end{small}
\end{thm}
The point is that, although it is not obvious how an RKHS structure arises from time-dependent kernels simultaneously over $t$, still a representer theorem applies. The consequence is that the optimization task is feasible over the infinite set $\cG$, i.e.~it is sufficient to minimize over a list of $T\times N$ parameters $a_{i,k}$. The proof is presented in the appendix.
\begin{rem}
In case the loss function is
$$loss(f(x),y)=(f(x)-y)^2,$$
the kernel regression problem $\min_{f\in\cG_{rep}}R_{emp}^\mathcal{\cS}(f)$ can be solved explicitly in terms of the Moore-Penrose pseudo inverse. Let $\vec{a},\ \vec{y}\in\mathbb{R}^{N\times T}$ denote the vectors of entries $a_{i,k}$ and $y_i^{(k)}$ and let ${\cK}\in Mat(N T\times N T,\mathbb{R})$ be the matrix whose entries are $\cK_{t_k}(x_i^{(k)},x_j^{(l)})$. The system of linear equations
$$\vec{y}={\cK}\cdot\vec{a},$$
can have no, a unique, or an infinite number of solutions. In either case the best match in terms of smallest $l^2$-norm is given by
$$\vec{a}^*={\cK}^+\cdot\vec{y},$$
where $\cK^+$ denotes the Moore-Penrose pseudo-inverse of $\cK$.

%Recall that in practice is it usually advisable to work with $QR$-decomposition instead of $\cK^+$ directly.
\end{rem}
\begin{example}(Minimal example)\label{minimalExample}
Suppose measurements outcomes 
\begin{align*}
\cS_1&=\{(x_{1}^{(1)}=1,y_{1}^{(1)}=1),(x_{2}^{(1)}=2,y_{2}^{(1)}=1)\},\\
\cS_2&=\{(x_{1}^{(2)}=1,y_{1}^{(2)}=1),(x_{2}^{(2)}=2,y_{2}^{(2)}=2)\}
\end{align*}
are obtained at times $t_1=1$ and $t_2=2$ from a function of the form $f(x,t)=\sum a_\mu \cK_t^{Gauss}(x,x_\mu)$. Suppose the 2-norm is used to measure loss. We have the following representation for the matrix $\cK= \cK^{Gauss}$
\begin{footnotesize}
\begin{align*}
&\left(\begin{array}{cccc}
\cK_{t_1}(x_1^{(1)},x_1^{(1)}) & \cK_{t_1}(x_1^{(1)},x_2^{(1)}) & ...& \cK_{t_1}(x_1^{(1)},x_2^{(2)})\\
\cK_{t_1}(x_2^{(1)},x_1^{(1)}) & \cK_{t_1}(x_2^{(1)},x_2^{(1)}) & ...& \cK_{t_1}(x_2^{(1)},x_2^{(2)})\\
\cK_{t_2}(x_1^{(2)},x_1^{(1)}) & \cK_{t_2}(x_1^{(2)},x_2^{(1)}) & ...& \cK_{t_2}(x_1^{(2)},x_2^{(2)})\\
\cK_{t_2}(x_2^{(2)},x_1^{(1)}) & \cK_{t_2}(x_2^{(2)},x_2^{(1)}) & ...& 
\cK_{t_2}(x_2^{(2)},x_2^{(2)})
\end{array}\right)
\end{align*}
\end{footnotesize}
where
$\vec{y}=\left(
y_{1}^{(1)}\ y_{2}^{(1)}\ y_{1}^{(2)}\ y_{2}^{(2)}\right)^T$
and 
$\vec{a}=\Big(
a_{1,1}\ a_{1,2}\ a_{2,1}\ a_{2,2}\Big)^T$. Computing $(\cK^{Gauss})^+\cdot\vec{y}$ using appropriate software yields the optimal coefficients $\vec{a}\approx\Big(
0.505\ 1.5984\ 0.505\ 1.5984\Big)^T$.
\end{example}
The theorem assumes a representation of candidate functions $f(x,t)$ in terms of infinite series (corresponding to a countable set $\mathbb{X}$). The convolution integral representation~\eqref{convIntegral} of PDE solutions is a continuous version of~\eqref{main}. The techniques for formalizing this generalization are, of course, well established. \cite[Section 5]{Christmann2008} contains a rigorous discussion.
To apply Thm.~\ref{thmRep} to the Fokker-Planck equation we represent its solution in terms of the convolution integral.
\begin{thm}\label{repThmDEq}
Let $\cF$ be the set of solutions of the Fokker-Planck equation with
$$L=-\partial_x\mu+\partial_{xx}D,$$
$D=D(x)$ and $\mu=\mu(x)$ and boundary conditions such that $L$ satisfies assumption~\emph{(A)}. Suppose that at times $t_k>0$, $k\in\{1,...,T\}$, data samples $\cS^{(k)}=\{(x_i^{(k)},y_i^{(k)})\}_{i=1,...,N}$ are measured.
Then
\begin{align*}
\inf_{f\in\cF}R_{emp}^\mathcal{\cS}(f)=\min_{f\in\cF_{rep}}R_{emp}^\mathcal{\cS}(f),
\end{align*}
where
\begin{small}
$$\cF_{rep}=\left\{f=f(x,t) |f(x,t)=\sum_{k=1}^T\sum_{i=1}^N a_{i,k} \cK^{FP}_t(x,x_i^{(k)})\right\}.$$
\end{small}
\end{thm}
\begin{figure}
%\begin{figure}[H]
\centering
\begin{subfigure}[b]{0.25\textwidth}
\includegraphics[width=\linewidth]{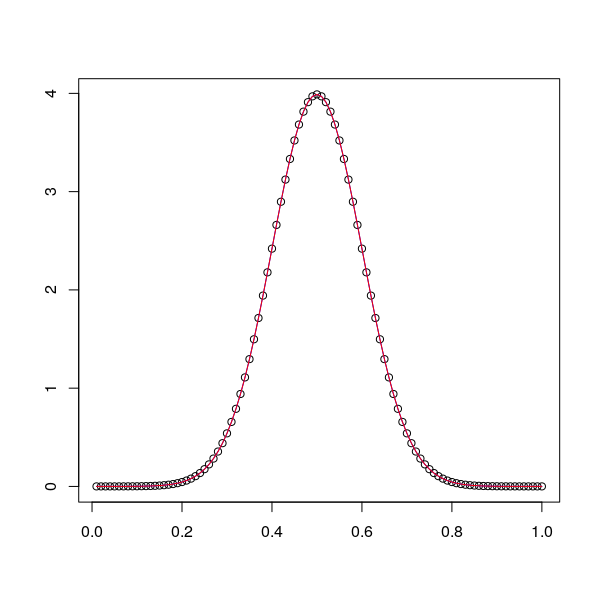}
\caption{ $t_0=0$}
\end{subfigure}
\begin{subfigure}[b]{0.25\textwidth}
\includegraphics[width=\linewidth]{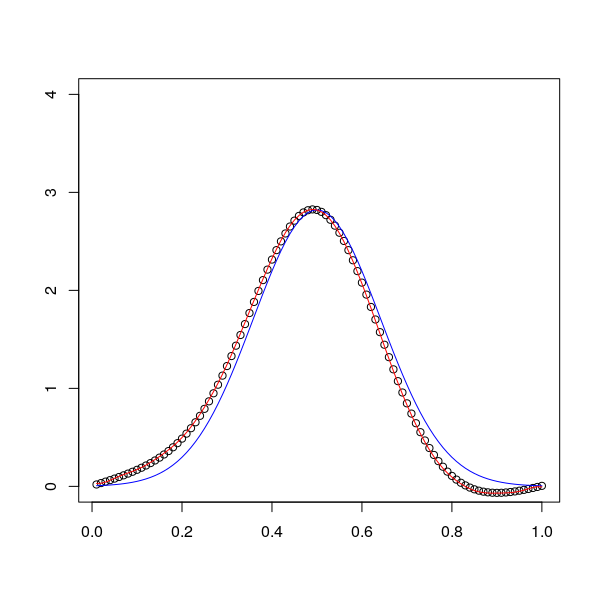}
\caption{$t_1=0.01$}
\end{subfigure}
\begin{subfigure}[b]{0.25\textwidth}
\includegraphics[width=\linewidth]{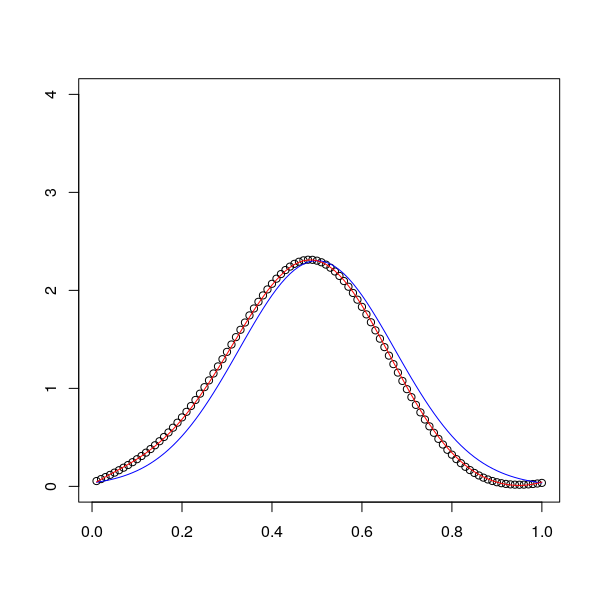}
\caption{$t_2=0.02$}
\end{subfigure}
\caption{Optimal match to temperature measurements at $100$ sensors and differnt times $t_0, t_1, t_2$ with initial condition $g(x)=1/{\sqrt{2\pi}}e^{-(x-1/2)^2/2}$. Black circles represent exemplary measurement outcomes. Blue lines depict temperature graphs as obtained by evolving the initial condition according to the heat equation. Red lines depict graphs of a regression function that exactly solves the heat equation.}\label{exampleHeatInit}
\end{figure}
The theorem links two different optimization exercises related to the PDE. On the one hand, there is the original optimization, which asks to minimize the empirical risk~\emph{over solutions of the Fokker-Planck equation}. On the other hand, there is the optimization of over functions $f(x,t)\in\cF_{rep}$
i.e.~over a~\emph{set of coefficients} $\{a_{i,k}\}_{i=1,...,N, k=1,...,T}$. By the Green's function property of the kernels, the coefficients reflect the initial conditions $f(x,t=0)$. In other words the optimization in Thm.~\ref{repThmDEq} is, de facto, an optimization over initial conditions. If the PDE is given with a priori initial conditions of the form $\{f(x,t=0)=g(x),\: x\in\mathbb{X}\}$ then its solution is unique and the optimization is, in principle, trivial. However, in practice initial conditions can be furnished with small errors (e.g.~from a measurement process). Such initial conditions can be naturally incorporated into the setting of Thm.~\ref{repThmDEq} by adding a data sample 
of the form $\{f(x_i^{(0)},t_0=0)=y_i^{(0)}\}_{i=1,...,N}$ to $\cS$. In this case the empirical risk decomposes as
\begin{small}
\begin{align}
&\sum_{i=1}^N loss\left(g(x_i^{(0)}),y_i^{(0)}\right)+\sum_{k=1}^T\sum_{i=1}^N loss\left(f(x_i^{(k)},t_k),y_i^{(k)}\right)\label{withInitConditions}
\end{align}
\end{small}
and Thm.~\ref{repThmDEq} is applied to the second term.
%
%
%
%
%Attempting to find the best regression function is equivalent to finding the best match initial conditions. The developed framework generalizes this observation in that it equally applies to over-determined problems.
%
%
%
\begin{example}(Heat Equation with Initial Conditions)\label{exHeatInit}
Suppose the temperature of a metal rod is measured are times $t_1=0.01$, $t_2=0.02$. The initial condition $\{f(x,t=0)=\frac{1}{\sqrt{2\pi}}e^{-(x-1/2)^2/2}\}$ is interpreted as a soft condition (and represented as a measurement at $t_0=0$). $100$ temperature sensors are placed at equal distances over the interval $[0,1]$. It is assumed that temperature evolves according to the heat equation and the $2$-norm measures loss. For illustration we suppose a measurement error of the form $0.2*\sin(2\pi x)$ at $t_1$ and $t_2$. Measurement outcomes and the optimal match solution obtained by minimizing~\eqref{withInitConditions} are shown in Figure~\ref{exampleHeatInit}.
\end{example}

Boundary conditions can be treated in similar vein. Hard boundary conditions imply that the kernels $\cK_t$ comply with the boundary conditions exactly. If a certain level of error can be tolerated, Thm.~\ref{repThmDEq} can be used to trade off between regression error and boundary conditions.
\begin{example}(Heat Equation with hard Dirichlet Boundary Conditions)\label{Ex:exampleHeatBoundHard}
Suppose, as before, that $100$ temperature sensors are equally placed over a metal rod and suppose measurements occur at times $t_1=0.01$, $t_2=0.02$, $t_3=0.03$. Suppose the rod is restricted to an interval $[0,1]$ and has temperature $0$ at $\{0,1\}$. As before temperature evolves according to the heat equation and the $2$-norm measures loss. For illustration it is assumed that measurement outcomes can be described by the functions $0.5-\abs{x-0.5}$ at $t_1$, by $0.3-0.6*\abs{x-0.5}$ at $t_2$ and by $0.2-0.4*\abs{x-0.5}$ at $t_3$, see Figure~\ref{introExample}. We apply Thm.~\ref{repThmDEq} to identify the optimal match solution to the heat equation with Dirichlet conditions. Straight forward computations show that
$$\varphi_n=\sqrt{2}\sin(n\pi x),\ \lambda_n=n^2\pi^2,$$
$$\cK_t^{DirichBound}(x,x')=2\sum_{n\geq0}e^{-n^2\pi^2t}\sin(n\pi x)\sin(n\pi x').$$
Figure~\ref{introExample} shows the optimal match solution obtained from Thm.~\ref{repThmDEq} using the pseudo-inverse.
%
%
%\begin{figure}
%\begin{figure}[H]
%\centering
%\begin{subfigure}[b]{0.2\textwidth}
%\includegraphics[width=\linewidth]{example3_0.png}
%\caption{$t_1=0.01$}
%\end{subfigure}
%\begin{subfigure}[b]{0.2\textwidth}
%\includegraphics[width=\linewidth]{example3_1.png}
%\caption{$t_2=0.02$}
%\end{subfigure}
%\caption{Optimal match to temperature measurements at $100$ sensors and differnt times $t_1, t_2$. Black circles represent exemplary measurement outcomes. Red lines depict graphs of an optimal match function that exactly solves the heat equation and complies with Dirichlet boundary.}\label{exampleHeatBoundHard}
%\end{figure}
\end{example}
In the real world dynamics the introduced framework can be expected to perform better in terms of~\emph{prediction} than agnostic techniques like IDW. This is illustrated by the following example. 
\begin{example}(Prediction under hard Dirichlet Boundary Conditions)
Suppose that as in Example~\ref{Ex:exampleHeatBoundHard} temperature measurements are taken at times $t_1=0.01$, $t_2=0.02$, $t_3=0.03$ from a metal rod with hard Dirichlet boundary conditions. Our goal is to predict the temperature at the future point in time $t_4=0.06$ given high levels of measurement error in the individual samples. The latter is modeled by the function $0.2*\sin(2\pi x)$ and added to each measurement. Two predictors are set up for comparison. The first is computed from Thm.~8. The second is the ordinary Gaussian kernel predictor $\cK^{Gauss}_s((x_i^{(k)},t_k),(x_j^{(l)},t_l))$, where the time-coordinates have been included as part of the observation. The bandwidth $s$ has been calibrated via cross validation to $0.45$. Figure~\ref{introExample} shows graphs of matches to data and respective predictions. The PDE predictor averages measurement error at different times to produce an accurate estimate of the exact evolution. 
\label{exampleIDW}
\end{example}
%
%
%
%\begin{example}(Heat Equation with soft Dirichlet Boundary Conditions)
%This example is identical to Example~\ref{Ex:exampleHeatBoundHard} above but boundary conditions need only be satisfied up to a measurement error.
%
%
%Figure~\ref{exampleHeatBoundSoft} shows the optimal match solution obtained.
%
%
%\begin{figure}
%\begin{figure}[H]
%\centering
%\begin{subfigure}[b]{0.2\textwidth}
%\includegraphics[width=\linewidth]{example4_0.png}
%\caption{$t_1=0.01$}
%\end{subfigure}
%\begin{subfigure}[b]{0.2\textwidth}
%\includegraphics[width=\linewidth]{example4_1.png}
%\caption{$t_2=0.02$}
%\end{subfigure}
%\caption{Optimal match to temperature measurements with soft boundary conditions, compare Figure~\ref{exampleHeatBoundHard}.}
%\label{exampleHeatBoundSoft}
%\end{figure}
%\end{example}

\subsection{Learning a time-dependent density}\label{densityLearning}
Suppose that a time-dependent density $\rho=\rho(t,x)$ evolves according to~\eqref{FokkerPlanck} and that samples $\cS^{(k)}=\{x_i^{(k)}\}_{i=1,...,N}$ are taken from $\rho_{t_k}=\rho(t_k,\cdot)$ at times $t_k$, $k\in\{1,...,T\}$. We construct a kernel density estimator that simultaneously takes account of all samples $\cS=\cup_{k=1}^T\cS^{(k)}$ while retaining consistency with~\eqref{FokkerPlanck}. To this aim we embed the initial density $\rho=\rho(x,t=0)$ into the RKHS of the Fokker-Planck kernels $\cK_t^{FP}$, representing it as
\begin{align*}
[\varepsilon(\rho)](x,t)=\int_{\mathbb{X}} \cK_t^{FP}(x,x')\rho(x')\d x'\in\cH_{\cK_t}.
\end{align*}
Notice that this convolution integral automatically solves the Fokker-Planck equation. Thus the embedding lifts the evolution of the density to the level of the RKHS and the developed regression theory applies. We generalize the KME by choosing the minimizer of
\begin{align}\label{densityLearningMinimization}
\cR_{KME}^{\cS}(\varepsilon)=\frac{1}{T}\sum_{k=1}^T\cR_{KME}^{\cS^{(k)}}(\varepsilon(\cdot,t_k))
\end{align}
as a \lq\lq{}simultaneous density estimator\rq\rq{}.
This functional arises from the time-dependent empirical risk~\eqref{regressionLearning} in the same way as the KME risk functional~\eqref{KME} arises from the ordinary (time-independent) empirical risk~\eqref{empiricalRisk}.
%
%
%
%
%
%
%
%
%
%
%\begin{align*}
%\hat{\rho}(x)=\frac{1}{N}\sum_{i=1}^N \cK(x,X_i).
%\end{align*}
%
%
%
%
%
%
%
%\begin{align*}
%\mathbb{E}\hat{\rho}(x)=\mathbb{E} \cK(x,X)=\int_\mathbb{X}\cK(x,y)\rho(y)\d y=\mu(\rho)(x).
%\end{align*}
%
%
%
%
%
%
%\begin{align*}
%\hat\rho(x)=\hat{\mu}(X_1,...,X_N)
%\end{align*}
%
%
%
%$$MSE=\mathbb{E}(\hat{\rho}(x)-\rho(x))^2$$
%
%
%Letzteres ist der Kerndicht Schaetzer.
%Nun was ist MISE?
%1) das ist die typische groesse die zu minimieren %es gilt fuer eine Dichte schaetzung.
%2) 
%$$MISE=\mathbb{E}\left(\int(\hat{\rho}(x)-\rho(x))^2\right)$$
%
%3) was machen den botev und so mit dem MISE?
%was ist was? $X$ gegeben mit Dichte? alles was man hat fuer MISE ist doch nur ein sample!!--> genau uns selbst das braucht man nicht. ein sample ist gar nicht noetig...
%
%
%
%
%
%
%
%
%
%
%
%
\begin{thm}[Time-dependent representer theorem for density estimation] Let $\cP$ be the set of densities $\rho(x,t)$ that solve the Fokker-Planck equation with
$$L=-\partial_x\mu+\partial_{xx}D,$$
$D=D(x)$ and $\mu=\mu(x)$ and boundary conditions such that $L$ satisfies assumption~\emph{(A)}.
Suppose that at times $t_k$, $k\in\{1,...,T\}$, samples $\cS^{(k)}=\{x_i^{(k)}\}_{i=1,...,N}$ are taken from the density $\rho(t_k)$. Then
\begin{align*}
\inf_{\rho\in\cP}R_{KME}^\mathcal{\cS}(\varepsilon(\rho))=\min_{\varepsilon\in\cP_{rep}}R_{KME}^\mathcal{\cS}(\varepsilon),
\end{align*}
where
\begin{small}
$$\cP_{rep}=\left\{\varepsilon=\varepsilon(x,t)| \varepsilon(x,t)=\sum_{k=1}^T\sum_{i=1}^N \beta_{i,k} \cK_t^{FP}(x,x_i^{(k)})\right\}.$$\label{repFullDensityLearning}
\end{small}
\end{thm}
Thm.~\ref{repFullDensityLearning} is an immediate consequence of Thm.~\ref{repThmDEq}, which is applied to~\eqref{densityLearningMinimization} to demonstrate that the minimizer is as a finite sum of kernels.
\begin{rem} Similar to regression, the optimization $\min_{\varepsilon\in\cP_{rep}}R_{KME}^\mathcal{\cS}(\varepsilon)$ can be solved explicitly in terms of the pseudo-inverse. First suppose $T=1$, let ${\cK}\in Mat(N\times N,\mathbb{R})$ be the matrix whose entries are $\cK_{ij}=\cK(x_i,x_j)$ and let $\cL=\sqrt{\cK}$. Let $\vec{\beta},\ \vec{1/N}\in\mathbb{R}^{N}$ denote vectors of entries $\beta_{i}$ and $1/N$. A quick computation shows that
\begin{align*}
\min_{\varepsilon\in\cP_{rep}}R_{KME}^{\cS}(\varepsilon)&=\min_{\vec{\beta}}(\vec{\beta}^T\cK\vec{\beta}-2*\vec{1/N}^T\cK\vec{\beta})\\
&=\min_{\vec{\beta}}\norm{\cL\vec{1/N}-\cL\vec{\beta}}{2}^2.
\end{align*}
Of course the minimum is $\vec{\beta}^*= \cL^+\cL\vec{1/N}=\vec{1/N}$ as mentioned before. This reasoning applies mutatis mutandis when $T>1$: The minimum of
$$\min_{\vec{\beta}}\sum_{k=1}^T\norm{\cL_{t_k}\vec{1/N}_k-\cL_{t_k}\vec{\beta}}{2}^2,$$
where now ${\cK}_{t}\in Mat(NT\times NT,\mathbb{R})$, $({\cK}_{t})_{ikjl}=\cK_t(x_i^{(k)},x_j^{(l)})$, ${\cL}_{t}=\sqrt{{\cK}_{t}}$ and $\vec{1/N}_k$ has entries $1/N$ corresponding to $\cS^{(k)}$ and $0$ else, can be computed using the pseudo-inverse of $(\cL_{t_1}\ \cL_{t_2}\ ...\ \cL_{t_T})^T$. 
\end{rem}
As before the optimization over all elements of $\cP$ is an optimization over possible choices of initial density. The optimal solution
$\varepsilon^*(x,t)$ corresponds to a list of coefficients $\{\beta_{i,k}^*\}$ and reflects the initial conditions 
$\rho(x,t=0)=\sum_{i,k}\beta_{i,k}^*\delta(x-x_i^{(k)})$.
Given these initial conditions the solution of the PDE is unique.

\begin{figure}
%\begin{figure}[H]
\centering
\begin{subfigure}[b]{0.25\textwidth}
\includegraphics[width=\linewidth]{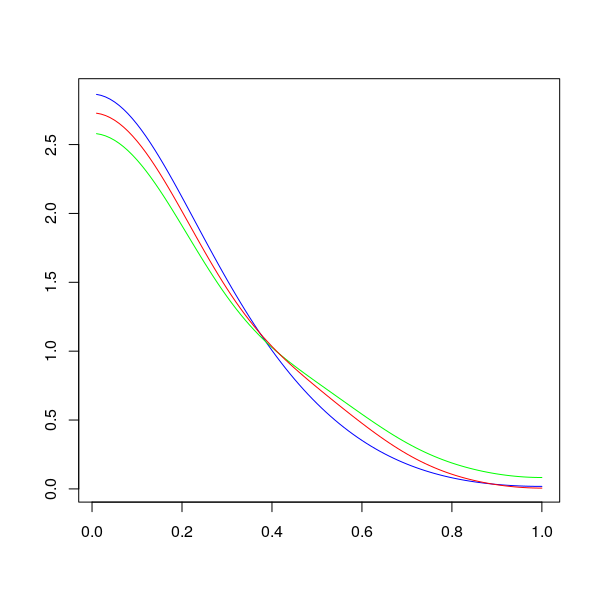}
\caption{ $t_1=0.01$}
\end{subfigure}
\begin{subfigure}[b]{0.25\textwidth}
\includegraphics[width=\linewidth]{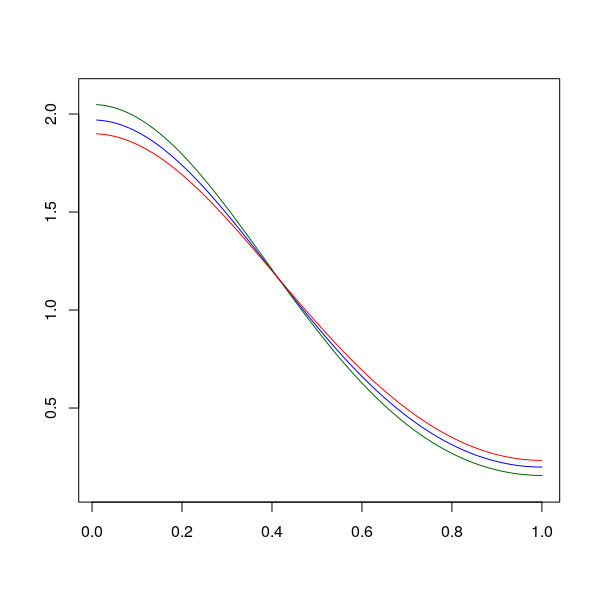}
\caption{$t_2=0.05$}
\end{subfigure}
\begin{subfigure}[b]{0.25\textwidth}
\includegraphics[width=\linewidth]{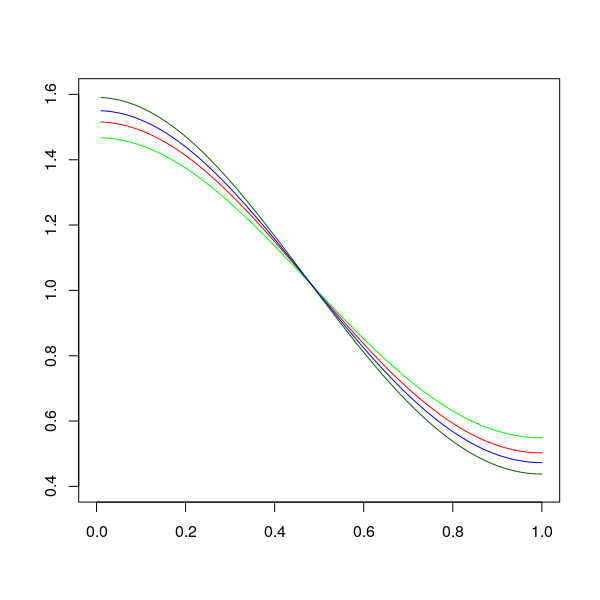}
\caption{$t_3=0.1$}
\end{subfigure}
\caption{Strenghening of individual kernel density estimators. Blue line deptics the evolution of the beta density. Bright green line depitcs $\hat{\rho}^{(NeumBound)}(x,t)$ constructed from $\cS^{(1)}$, dark green line depicts $\hat{\rho}^{(NeumBound)}(x,t)$ from $\cS^{(2)}$. Red line depitcs the combined estimator from $\cS$.}\label{fig:density}
\end{figure}

A common situation in kernel density estimation is that the domain $\mathbb{X}$ of data is known in advance, say~$\mathbb{X}=[0,1]$. This leads to the topic of boundary conditions. The homogeneous Neumann boundary conditions
\begin{align*}
\partial_x\rho(x,t)|_{x=0}=\partial_x\rho(x,t)|_{x=1}=0
\end{align*}
ensure that $\partial_t\int_{\mathbb{X}}\rho(x,t)\d x =0$, which entails 
$$\int_{\mathbb{X}}\rho(x,t)\d x =\int_{\mathbb{X}}\rho(x,0)\d x =1.$$ 
\begin{example}(Learning an evolving density on $[0,1]$)
Suppose that two samples each of size $N=100$ are taken at times $t_1=0.01$ and $t_2=0.05$ from a density that evolves according to the heat equation with homogeneous Neumann boundary conditions on $[0,1]$. For illustration we suppose that at time $t=0$ the density is the beta density $\rho(x,t=0)=4(1-x)^3$ and we construct an estimator for the density at $t_3=0.1$ that accounts for both samples.

The kernel of the heat equation under Neumann boundary conditions is
\small
\begin{align*}
\cK^{NeumBound}_t(x,x')=1+2\sum_{n>0}e^{-n^2\pi^2t}\cos(n\pi x)\cos(n\pi x').
\end{align*}
\normalsize
The individual estimators for samples at times $t_1$ and $t_2$ are $\hat{\rho}^{(NeumBound)}(x,t_1)$ and $\hat{\rho}^{(NeumBound)}(x,t_2)$, where $\hat\rho^{NeumBound}(x,t)=1/N\sum_{i=1}^N\cK^{NeumBound}_t(x,X_i)$. A combined estimator for the evolving density with coefficients $\{\beta_{i,1}^*,\beta_{i,2}^*\}_{i=1,...,100}$ is provided by Thm.~\ref{repFullDensityLearning}, see Figure~\ref{fig:density}.
\end{example}

\section{Conclusion}
Time-dependent PDEs occur in countless situations throughout scientific disciplines. Kernel methods are commonly employed in PDE theory~\cite{Schaback2006} but the existing techniques mostly focus on the time-independent context. Similarly, static kernel methods are common in spatio-temporal modeling and fall under the general IDW paradigm. Models that accurately take account of the system dynamics are much less common. The article at hand introduces a new kernel-based paradigm for spatio-temporal modeling based on time-dependent kernels that realize the Fokker-Planck dynamics. A respective kernel-based learning theory for time-dependent PDEs is introduced and a representer theorem is provided for the application of dynamic kernel techniques. 

Our kernel density estimators are closely related to the famous~{diffusion-estimator} of~\cite{Botev2010}. The latter uses the Fokker-Planck equation as a resource for the construction of estimators under prior information, such as domain constrains and a pilot density estimate. On the theoretical side we expect applications of our method in the construction of diffusion estimators from multiple correlated samples along the lines of~\cite{Botev2010}.

\section{Appendix}\label{Sec:appendix}
\begin{proof}[Proof of Lemma~\ref{positiveKernel}]
Symmetry of the kernel, i.e. $\cK^{FP}_t(x,x')=\cK^{FP}_t(x',x)$ follows directly by plugging in. To show positive-definiteness set $u = \sum_ic_i\varphi_n(x_i)$ and  compute
$
\sum_{ij}^Nc_ic_j\cK^{FP}_t(x_i,x_j)=
\sum_{n=1}^\infty e^{-\lambda_nt}uu^T,
$ which is positive-definite because $e^{-\lambda_nt}>0$.
\end{proof}

\begin{proof}[Proof of Thm.~\ref{thmRep}]
By definition $\cG_{rep}\subset\cG$, i.e.~
\begin{align*}
\inf_{f\in\cG}R_{emp}^\mathcal{\cS}(f)\leq\min_{f\in\cG_{rep}}R_{emp}^\mathcal{\cS}(f).
\end{align*}
Fix $t_k$ and let, by the Moore-Aronszajn theorem, $\cH_{t_k}$ be the unique reproducing kernel Hilbert space with kernel $\cK_{t_k}$. Evaluating $f\in\cG$ at $t_k$ simply gives
$
f(x,t_k)=\sum_{\mu=1}^\infty a_\mu \cK_{t_k}(x,x_\mu)$.
We view (for fixed $t_k$) the function $f(x,t_k)$ as an element of $\cH_{t_k}$. To emphasize this we introduce the notation $f_{t_k}\in\cH_{t_k}$, i.e.~$f(x,t_k)=f_{t_k}(x)$. Let $P_{t_k}$ be the projector in $\cH_{t_k}$ on 
\begin{align*}
K=span\left\{\cK_{t_k}(x,x_i^{(k)})\right\}_{i=1,...,N}
\end{align*}
and let $Q_{t_k}$ be the projector onto the orthogonal complement $P_{t_k}+Q_{t_k}={1}$. For any function $g\in\cH_{t_k}$, which is contained in the orthogonal complement of $K$ the reproducing property of $\cK_{t_k}$ gives
$$0=\braket{g(x)}{\cK_{t_k}(x,x_i^{(k)})}_{Aronszajn}=g(x_i^{(k)}).$$
This implies $(Q_{t_k}f)(x_i^{(k)},t_k)=0.$ As a consequence we find that
\begin{align}
f(x_i^{(k)},t_k)=f_{t_k}(x_i^{(k)})=(P_{t_k}f_{t_k})(x_i^{(k)})+(Q_{t_k}f_{t_k})(x_i^{(k)})=(P_{t_k}f_{t_k})(x_i^{(k)})\label{proj}.
\end{align}
Notice that 
\begin{align}
(P_{t_k}f_{t_k})(x)=\sum_{j=1}^Na_{j}^{(k)} \cK_{t_k}(x,x_j^{(k)}),\label{projected}
\end{align}
where $\{a_{j}^{(k)}\}_{j=1,...,N}$ constitute a subset of the coefficients $\{a_\mu\}_{\mu=1,...,\infty}$ corresponding to the basis vectors $\{\cK_{t_k}(x,x_j^{(k)})\}_{j=1,...,N}$. The coefficients $a_{j}^{(k)}$ do not depend on $t$ since, by assumption, the coefficients $a_\mu$ do not depend on $t$. Taking together equations~\eqref{proj} and~\eqref{projected} for any fixed $t_k$ it follows that
\begin{small}
\begin{align*}
&loss\left(f(x_i^{(k)},t_k),y_i^{(k)}\right)
= loss\left(\sum_{j=1}^Na_{j}^{(k)} \cK_{t_k}(x_i^{(k)},x_j^{(k)}),y_i^{(k)}\right).
\end{align*}
\end{small}
Since the above holds for any $t_k$ we obtain
\begin{align}
\inf_{f\in\cG}R_{emp}^\mathcal{\cS}(f)=
\min_{\cup_{j,k}\{a_{j}^{(k)}\}}\sum_{k=1}^T\sum_{i=1}^N loss\left(\sum_{j=1}^Na_{j}^{(k)} \cK_{t_k}(x_i^{(k)},x_j^{(k)}),y_i^{(k)}\right).\label{infProb}
\end{align}
On the other hand by assumption any function in $\cG_{rep}$ is of the form
$$f(x,t)=\sum_{k=1}^T\sum_{j=1}^N c_{j,k} \cK_{t}(x,x_j^{(k)}),$$
where we wrote $c_{j,k}$ (instead of $a_{j,k}$) to distinguish them from $a_{j}^{(k)}$. Plugging in we find
\begin{align}
\min_{f\in\cG_{rep}}R_{emp}^\mathcal{\cS}(f)=\min_{\cup_{j,l}\{c_{j,l}\}}\sum_{k=1}^T\sum_{i=1}^N loss\left(\sum_{l=1}^T\sum_{i=1}^Nc_{j,l} \cK_{t_k}(x_i^{(k)},x_j^{(l)}),y_i^{(k)}\right).\label{minProb}
\end{align}
Comparing~\eqref{infProb} and~\eqref{minProb}
it follows that
\begin{align*}
\inf_{f\in\cG}R_{emp}^\mathcal{\cS}(f)\geq\min_{f\in\cG_{rep}}R_{emp}^\mathcal{\cS}(f)
\end{align*}
because any candidate function for the left hand side is also a candidate function for the right hand side.

\end{proof}
\begin{proof}[Proof of Thm.~\ref{repThmDEq}]
Following the discussion of the Fokker-Planck equation any solution $f\in\cF$ can be written in the form
$$f(x,t)=\int_{\mathbb{X}}g(x')K_t^{FP}(x,x')\d x'.$$
By Lem.~\ref{positiveKernel}, $K_t^{FP}(x,x')$ is a symmetric, positive-definite kernel for $t>0$. As in the proof of Thm.~\ref{thmRep} it generates respective RKHS $\cH_{t_k}$ at times $t_k$. Let now $f(x,t)\in\cF$. To arrive at the representation~\eqref{projected} we apply Thm.~5.5 of~\cite{Christmann2008}. The proof then follows the same line as the proof of Thm.~\ref{thmRep}.
\end{proof}

\bibliographystyle{plain}
\bibliography{bibfile1}

\end{document}